\documentclass{article}

\usepackage{amsfonts,amssymb,amsmath,rotating}
\allowdisplaybreaks[3]

\newcommand{\ncd}{\newcommand}
\ncd{\gvct}[1]{{\boldsymbol{#1}}}
\def\N{\mathbb N}
\def\R{\mathbb R}
\def\Z{\mathbb Z}
\def\BAL#1\EAL{\begin{align}#1\end{align}}
\ncd{\ip}[2]{\langle #1 | #2\rangle}
\ncd{\Eq}[1]{Eq.~(\ref{#1})}

\ncd{\CG}{{\cal G}}
\ncd{\CH}{{\cal H}}
\DeclareMathOperator{\den}{den}
\DeclareMathOperator{\lcm}{lcm}
\newtheorem{theo}{Theorem}[section]

\begin{document}

\title{Coincidences of Hypercubic Lattices in 4 dimensions}

\author{P. Zeiner\\
Institute for Theoretical Physics \& CMS, TU Wien,\\
Wiedner Hauptsra{\ss}e 8--10, 1040 Vienna, Austria}

\maketitle

\begin{abstract}
We consider the CSLs of $4$--dimensional hypercubic lattices.
In particular,
we derive the coincidence index $\Sigma$ and calculate the number of
different CSLs as well as the number of inequivalent CSLs for a given $\Sigma$.
The hypercubic face centered case is dealt with in detail and it is sketched
how to derive the corresponding results for the primitive hypercubic lattice.
\end{abstract}

\section{Introduction}

Coincidence site lattices (CSL) for three--dimensional lattices
have been studied intensively since they are an important tool to characterize
and
analyze the structure of grain boundaries in crystals~(\cite{boll70,boll82}
and references therein). For quasicrystals these concepts have to be adapted.
Since a lot of quasiperiodic structures can be obtained by the well--known
cut and projection scheme~\cite{dun1,mbaake98}
from a periodic structure in superspace, it is
natural to also investigate CSLs in higher dimension. An important example
are the four--dimensional hypercubic lattices, which shall be discussed here.
Four--dimensional lattices are particularly interesting since they are the
first ones that allow 5--fold, 8--fold, 10--fold and 12--fold symmetries
which are actually observed in quasicrystals. In particular the
four--dimensional hypercubic lattices allow 8--fold symmetries, from which we
can obtain e.g. the prominent Ammann--Beenker tiling~\cite{mbaake98}.

Since rotations in four--dimensional space can be parameterized by quaternions,
one has a strong tool to investigate the CSLs of the hypercubic lattices.
In particular, one knows all coincidence rotations~\cite{baa97} and thus
all CSLs can be characterized. But one can go much further and this will be
done in the present paper. We first calculate the coincidence index $\Sigma$
and then we try to find all CSLs for a given index $\Sigma$. It turns out
that one can calculate the total number of different CSLs for a given
$\Sigma$ and furthermore, we can derive even the number of inequivalent
CSLs and for each CSL we can calculate the number of equivalent CSLs.

These calculations are facilitated by the fact that the four--dimensional
rotations are closely related to the their three--dimensional counterparts.
In particular we exploit the fact that $SO(4)\simeq SU(2)\times SU(2)/C_2$,
i.e. we can make use of the results of the three--dimensional cubic
case that have been published recently~\cite{pzcsl1}. Thus the results of
the hypercubic case are quite similar to the three--dimensional results,
although proofs are a bit more lengthy and the resulting formulas are a bit
more complex. However, there is one big difference between three and four
dimensions: Whereas all important quantities like $\Sigma$, number of CSLs
etc. are equal for all three kinds of cubic lattices, this is no longer true
for four dimensions. For a given coincidence rotation $R$, the coincidence
indices for the primitive and the face centered hypercubic lattice are in
general not the same, which is not surprising since the point groups are
different, too. Thus we must deal with both cases separately. However,
one can derive the results of the primitive lattice from the corresponding
results of the face centered lattice. Thus we concentrate on the latter and
sketch how these results can then be used for the primitive hypercubic
lattice.

Now let us recall some basic facts and fix the notation. 
Let $\gvct{L}\subseteq\R^n$ be an $n$-dimensional lattice and $R$ a rotation.
Then $\gvct{L}(R)=\gvct{L}\cap R\gvct{L}$ is called a
\emph{coincidence site
lattice}~(CSL) if it is a sublattice of finite index of $\gvct{L}$,
the corresponding rotation is called a coincidence rotation~\cite{baa97}. The
coincidence index $\Sigma(R)$ is defined as the index of $\gvct{L}(R)$ in
$\gvct{L}$. By index we mean the group theoretical index of $\gvct{L}(R)$ in
$\gvct{L}$, where we view $\gvct{L}(R)$ and $\gvct{L}$ as additive groups.

Any rotation in $4$ dimensions can be parameterized by two
quaternions $\gvct{p}=(k,\ell,m,n)$ and $\gvct{q}=(a,b,c,d)$ in the
following way~\cite{koecheng,hurw,val}:
\BAL
&R(\gvct{p},\gvct{q})=\frac{1}{|\gvct{p}\gvct{q}|}M(\gvct{p},\gvct{q})\\
\label{parrot}
&M(\gvct{p},\gvct{q})=
\begin{pmatrix}
\ip{\gvct{p}}{\gvct{q}}&
\ip{\gvct{p}\gvct{u}_1}{\gvct{q}}&
\ip{\gvct{p}\gvct{u}_2}{\gvct{q}}&
\ip{\gvct{p}\gvct{u}_3}{\gvct{q}}\\
\ip{\gvct{p}}{\gvct{u}_1\gvct{q}}&
\ip{\gvct{p}\gvct{u}_1}{\gvct{u}_1\gvct{q}}&
\ip{\gvct{p}\gvct{u}_2}{\gvct{u}_1\gvct{q}}&
\ip{\gvct{p}\gvct{u}_3}{\gvct{u}_1\gvct{q}}\\
\ip{\gvct{p}}{\gvct{u}_2\gvct{q}}&
\ip{\gvct{p}\gvct{u}_1}{\gvct{u}_2\gvct{q}}&
\ip{\gvct{p}\gvct{u}_2}{\gvct{u}_2\gvct{q}}&
\ip{\gvct{p}\gvct{u}_3}{\gvct{u}_2\gvct{q}}\\
\ip{\gvct{p}}{\gvct{u}_3\gvct{q}}&
\ip{\gvct{p}\gvct{u}_1}{\gvct{u}_3\gvct{q}}&
\ip{\gvct{p}\gvct{u}_2}{\gvct{u}_3\gvct{q}}&
\ip{\gvct{p}\gvct{u}_3}{\gvct{u}_3\gvct{q}}
\end{pmatrix}\\
&=
\begin{pmatrix}
ak+b\ell+cm+dn & -a\ell+bk+cn-dm & -am-bn+ck+d\ell & -an+bm-c\ell+dk\\
a\ell-bk+cn-dm & ak+b\ell-cm-dn & -an+bm+c\ell-dk & am+bn+ck+d\ell \\
am-bn-ck+d\ell & an+bm+c\ell+dk & ak-b\ell+cm-dn & -a\ell-bk+cn+dm \\
an+bm-c\ell-dk & -am+bn-ck+d\ell & a\ell+bk+cn+dm & ak-b\ell-cm+dn
\end{pmatrix}
\EAL
Here, $\gvct{u}_i$ are the unit quaternions $\gvct{u}_0=(1,0,0,0)$,
$\gvct{u}_1=(0,1,0,0)$, $\gvct{u}_2=(0,0,1,0)$ and $\gvct{u}_3=(0,0,0,1)$ and
$|\gvct{p}|^2=k^2+\ell^2+m^2+n^2$ is the norm of $\gvct{p}$. Furthermore
we have made use of the inner product
$\ip{\gvct{p}}{\gvct{q}}:=ak+b\ell+cm+dn$.
If we identify the quaternions with the elements of $\Z^4$ (or $\R^4$)
in the obvious way,
then the action of $M(\gvct{p},\gvct{q})$ on a vector $\gvct{x}\in\Z^4$ can
be written as $M(\gvct{p},\gvct{q})\gvct{x}=\gvct{p}\gvct{x}\gvct{\bar q}$.
Here $\gvct{\bar q}=(a,-b,-c,-d)$ denotes the conjugate of $\gvct{q}$.
By an \emph{integral} quaternion we
mean a quaternion with integral coefficients. If the greatest common divisor
of all coefficients is $1$, we call the quaternion \emph{primitive}.
In the following,
all quaternions will be either primitive or normalized to unity. It will
always be clear from the context which convention has been chosen. 
Obviously $R(\gvct{p},\gvct{q})$
is a rational matrix if $\gvct{p}$ and $\gvct{q}$ are integral quaternions
such that $|\gvct{p}\gvct{q}|$ is an integer. 
In this case we call the
pair $(\gvct{p},\gvct{q})$ \emph{admissible}~\cite{baa97}. 
On the other hand any
rational orthogonal matrix $R$ can be parameterized by an admissible pair
of integral quaternions. Furthermore recall
$\gvct{p}^{-1}=\frac{1}{|\gvct{p}|^2}\gvct{\bar p}$.


\section{The CSLs and their $\Sigma$--values}

In $4$ dimensions there are only two different hypercubic lattices, namely
the primitive and the centered hypercubic lattices. They are equivalent to
$\gvct{L}_P=\Z^4$ and $\gvct{L}_F=D_4$, respectively. $D_4\subset \Z^4$ is of
index 2 and consists of all integer vectors $\gvct{n}$ with
$|\gvct{n}|^2$ even. It is known that $R$ is a coincidence rotation
of $\Z^4$ or $D_4$, if and only if all its entries are rational~\cite{baa97},
i.e. $R=R(\gvct{p},\gvct{q})$ for some admissible
pair of primitive quaternions $(\gvct{p},\gvct{q})$. In order to analyze
the CSLs it is often convenient to find some appropriate sublattices
of the CSLs. To this end we define the \emph{denominator}
\BAL
\den(R)=\gcd\{k\in\N | k R \gvct{L} \subseteq \gvct{L}(R)\},
\EAL
where $\gcd$ denotes the greatest common divisor. Since $\den(R)\cdot \gvct{L}$
is a sublattice of $\gvct{L}(R)$
it follows that
\BAL\label{dens}
\den(R)\leq \Sigma(R)\leq\den(R)^n,
\EAL
where $n=4$ is the dimension of $\gvct{L}$. In case of the primitive cubic
lattice this definition coincides with~\cite{baa97}
\BAL
\den_P(R)=\gcd\{k\in\N | k R \mbox{ integer matrix}\},
\EAL
whereas for the centered lattice we find
\BAL
\den_F(R)=2^{-\ell}\den_P(R),
\EAL
where $\ell=0,1$ is the maximal power such that $2^\ell$ divides $\den_P(R)$.
In particular we find for any admissible pair $\gvct{p},\gvct{q}$
\BAL
\den_F(R(\gvct{p},\gvct{q}))=2^{-\ell}|\gvct{p}\gvct{q}|,
\EAL
where $\ell=0,1,2$ is the maximal power such that $2^\ell$ divides
$|\gvct{p}\gvct{q}|$.

It follows from \Eq{dens} that $R$ is a symmetry
operation of $\gvct{L}$ if and only if $\den(R)=1$. Thus $R(\gvct{p},\gvct{q})$
is a symmetry operation of the centered lattice $D_4$ if and only if
$|\gvct{p}|^2=1,4$ and
$|\gvct{q}|^2=1,4$ or $|\gvct{p}|^2=|\gvct{q}|^2=2$. This gives the well known
576 pure symmetry rotations of $D_4$. Note that not all of them are
integer matrices, which reflects the fact that the symmetry group of $D_4$
is larger than those of $\Z^4$.
In fact only 192 rotations are integer matrices, namely
the pure symmetry rotations of $\Z^4$. These are the rotations corresponding
to the pairs $(\gvct{p},\gvct{q})$ such that $|\gvct{p}|^2=|\gvct{q}|^2=1$ or
$|\gvct{p}|^2=|\gvct{q}|^2=2$ with $\ip{\gvct{p}}{\gvct{q}}$ even or
$|\gvct{p}|^2=|\gvct{q}|^2=4$ with $\ip{\gvct{p}}{\gvct{q}}$ divisible by $4$.

We consider the face centered lattice first, formulating and proving a result
that was first stated (without proof) in \cite{baa97}, Eq.~(3.21).
\begin{theo}
Let $\gvct{p},\gvct{q}$ be an admissible pair of primitive integer
quaternions and let $\Sigma(\gvct{p})=2^{-\ell}|\gvct{p}|^2$,
where $\ell=0,1,2$ is the maximal power such that $2^\ell$ divides
$|\gvct{p}|^2$. Then, for the fcc--lattice,
the rotation $R(\gvct{p},\gvct{q})$ has coincidence index
\BAL
\Sigma_F(\gvct{p},\gvct{q}):=\Sigma_F(R(\gvct{p},\gvct{q}))=
\lcm(\Sigma(\gvct{p}),\Sigma(\gvct{q})).
\EAL
\end{theo}
\emph{Proof:}
Let us write $|\gvct{p}|^2=\alpha^2\gamma$, $|\gvct{q}|^2=\beta^2\gamma$, where
$\gamma=\gcd(|\gvct{p}|^2,|\gvct{q}|^2)$. Further let
$\gvct{p}^{(i)}=\gvct{p}\gvct{u}_i$ and
$\gvct{\bar q}^{(i)}=\gvct{\bar u}_i\gvct{\bar q}$. Then $\beta\gvct{p}^{(i)}$
and $\alpha\gvct{\bar q}^{(j)}$ are integer vectors with integer pre-images.
Thus they are in $\gvct{L}_P(R)$ and hence certainly  $2\beta\gvct{p}^{(i)}$
and $2\alpha\gvct{\bar q}^{(j)}$ are in $\gvct{L}_F(R)$.\footnote{Unless
$|\gvct{p}|^2$ and $|\gvct{q}|^2$ are both odd, even $\beta\gvct{p}^{(i)}$
and $\alpha\gvct{\bar q}^{(j)}$ are elements of $\gvct{L}_F(R)$. In any case
$\beta(\gvct{p}^{(i)}+\gvct{p}^{(j)})\in \gvct{L}_F(R)$ and
$\alpha(\gvct{\bar q}^{(k)}+\gvct{\bar q}^{(\ell)}) \in \gvct{L}_F(R)$.}
Thus if $i\ne j$ and $k\ne \ell$ the four vectors $2\beta\gvct{p}^{(i)},
2\beta\gvct{p}^{(j)},2\alpha\gvct{\bar q}^{(k)},2\alpha\gvct{\bar q}^{(\ell)}$
span a sublattice of $\gvct{L}_F(R)$. Now
\BAL
\det\left(\beta\gvct{p}^{(i)},\beta\gvct{p}^{(j)},
\alpha\gvct{\bar q}^{(k)},\alpha\gvct{\bar q}^{(\ell)}\right)=
\alpha^2\beta^2(\ip{\gvct{p}^{(i)}}{\gvct{\bar q}^{(k')}}
\ip{\gvct{p}^{(j)}}{\gvct{\bar q}^{(\ell')}}-
\ip{\gvct{p}^{(i)}}{\gvct{\bar q}^{(\ell')}}
\ip{\gvct{p}^{(j)}}{\gvct{\bar q}^{(k')}},
\EAL
where $k',\ell'$ are chosen such that $(k',\ell',k,\ell)$ is an even
permutation of $(0,1,2,3)$.
Hence we conclude that $\Sigma_F(R)$ divides $8\alpha^2\beta^2c$, where $c$
is the greatest common divisor of
$c^{ijk\ell}=\ip{\gvct{p}^{(i)}}{\gvct{\bar q}^{(k)}}
\ip{\gvct{p}^{(j)}}{\gvct{\bar q}^{(\ell)}}-
\ip{\gvct{p}^{(i)}}{\gvct{\bar q}^{(\ell)}}
\ip{\gvct{p}^{(j)}}{\gvct{\bar q}^{(k)}})$. Using the expansion
\BAL
|\gvct{p}|^2\gvct{a}=\sum_{i=0}^{3}\ip{\gvct{p}^{(i)}}{\gvct{a}}\gvct{p}^{(i)}
\EAL
we see that $c$ divides
\BAL
\sum_{i=0}^{3}c^{ijk\ell}\ip{\gvct{p}^{(i)}}{\gvct{a}}=
|\gvct{p}|^2(\ip{\gvct{a}}{\gvct{\bar q}^{(k)}}
\ip{\gvct{p}^{(j)}}{\gvct{\bar q}^{(\ell)}}-
\ip{\gvct{a}}{\gvct{\bar q}^{(\ell)}}
\ip{\gvct{p}^{(j)}}{\gvct{\bar q}^{(k)}})
\EAL
for any integer quaternion $\gvct{a}$. We now choose $\gvct{a}$ such that
$\ip{\gvct{a}}{\gvct{\bar q}^{(k)}}=0$, in particular if $k=0$ we choose
$\gvct{a}=q_\ell\gvct{u}_0+q_0\gvct{u}_\ell$ and
$\gvct{a}=\gvct{\bar q}^{(\ell)}-q_\ell\gvct{u}_0-q_0\gvct{u}_\ell$.
Hence $c$ must divide
$|\gvct{p}|^2(q_m^2+q_n^2)\ip{\gvct{p}^{(j)}}{\gvct{\bar q}}$.
Now $\gvct{\bar q}$
is primitive so that the greatest common divisor of all combinations
$q_m^2+q_n^2$ is at most $2$. Thus $c$ divides
$2|\gvct{p}|^2\ip{\gvct{p}^{(j)}}{\gvct{\bar q}}$. Similarly one proofs that
$c$ divides $2|\gvct{p}|^2\ip{\gvct{p}^{(j)}}{\gvct{\bar q}^{(k)}}$
for arbitrary
$k$, and hence $c$ must divide $8|\gvct{p}|^2$ because $\gvct{q}$ and
$\gvct{p}$ are both primitive. In the same way one shows that $c$ divides
$8|\gvct{q}|^2$. Thus $c$ divides $8\gamma$ and $\Sigma_F(R)$ divides
$64\alpha^2\beta^2\gamma$. But $\Sigma_F(R)$ divides $\den_F(R)^4$,
which is odd. So $\Sigma_F(R)$ divides
$\lcm(\Sigma(\gvct{p}),\Sigma(\gvct{q}))$.

It remains to show the converse statement,
$\lcm(\Sigma(\gvct{p}),\Sigma(\gvct{q}))\leq \Sigma_F(R)$.
To this end, we count the vectors $\gvct{y}\in\gvct{L}_F(R)$ contained in the
hypercube $H(2\beta\gvct{p}^{(i)})$ spanned by $2\beta\gvct{p}^{(i)}$.
If there are $n_F$ of them then
$\Sigma_F(R)=8\beta^4|\gvct{p}|^4/n_F=8\alpha^4|\gvct{q}|^4/n_F$.
Now $\gvct{L}_F(R)$ is a sublattice of
$\gvct{L}_P(R)$, so that $\Sigma_F(R)$ is a multiple of
$8\beta^4|\gvct{p}|^4/n_P$ if $n_P$ denotes the number of vectors
$\gvct{y}\in\gvct{L}_P(R)$ contained in the hypercube
$H(2\beta\gvct{p}^{(i)})$.
Now $n_P=16n'_P$, where $n'_P$ is the number of the vectors
$\gvct{y}\in\gvct{L}_P(R)$ contained in the smaller hypercube
$H(\beta\gvct{p}^{(i)})$ spanned by $\beta\gvct{p}^{(i)}$. Equivalently
we can count their pre-images
$\gvct{x}=R^{-1}\gvct{y}$ lying inside the hypercube
$H(\alpha\gvct{u}_i\gvct{q})$. In the following we identify
$H(\alpha\gvct{u}_i\gvct{q})$ with the factor group
$\gvct{L}_P/\gvct{L}_{\gvct{q}}$, where $\gvct{L}_{\gvct{q}}$ denotes the
$\Z$--span of the vectors $\alpha\gvct{u}_i\gvct{q}$. 

Observe that any vector $\gvct{x}$ of $H(\alpha\gvct{u}_i\gvct{q})$ can be
expressed as
\BAL
\gvct{x}=\frac{1}{|\gvct{q}|^2}
\sum_{i=0}^{3}\ip{\gvct{x}}{\gvct{u}_i\gvct{q}}\gvct{u}_i\gvct{q},
\EAL
such that $0\leq\ip{\gvct{x}}{\gvct{u}_i\gvct{q}}<\alpha|\gvct{q}|^2$. Now
$\gvct{x}$ is in $\gvct{L}_P(R)$ if its image
\BAL\label{condRx}
R(\gvct{p},\gvct{q})\gvct{x}=\frac{1}{|\gvct{p}\gvct{q}|}
\sum_{i=0}^{3}\gvct{p}^{(i)}\ip{\gvct{x}}{\gvct{u}_i\gvct{q}}
\EAL
is an integral vector. Since
\BAL\label{condbeta}
\ip{\gvct{p}^{(i)}}{R\gvct{x}}=
\frac{|\gvct{p}|}{|\gvct{q}|}\ip{\gvct{x}}{\gvct{u}_i\gvct{q}}=
\frac{\alpha}{\beta}\ip{\gvct{x}}{\gvct{u}_i\gvct{q}}
\EAL
all coefficients $\ip{\gvct{x}}{\gvct{u}_i\gvct{q}}$ must be divisible by
$\beta$. In order to determine the vectors that satisfy this condition
we first observe that there exists a vector $\gvct{x}$
such that $\ip{\gvct{x}}{\gvct{q}}=1$ since $\gvct{q}$ is primitive. Regarding
$\gvct{x}$ as an element of the abelian group $\gvct{L}_P/\gvct{L}_{\gvct{q}}$
we see that it has order $\alpha|\gvct{q}|^2$. Among all vectors
$\gvct{x}'$ with $\ip{\gvct{x}'}{\gvct{q}}=0$ there exists one of order
$\alpha|\gvct{q}|^2/2$ or $\alpha|\gvct{q}|^2$, depending on whether
$|\gvct{q}|^2$ is divisible by $4$ or not.\footnote{Consider the vectors
$q_\ell\gvct{u}_0-q_0\gvct{u}_\ell$ and
$\gvct{u}_\ell\gvct{q}+q_\ell\gvct{u}_0-q_0\gvct{u}_\ell$. Their orders
are multiples of $\alpha|\gvct{q}|^2/\gcd(\alpha|\gvct{q}|^2,q_\ell^2+q_m^2)$
and an appropriate combination thereof gives the desired vector $\gvct{x}'$.}
Hence $\gvct{x}$ and $\gvct{x}'$ generate a subgroup of
order $\alpha^2|\gvct{q}|^4$ or
$\alpha^2|\gvct{q}|^4/2$ of $\gvct{L}_P/\gvct{L}_{\gvct{q}}$.
Condition~(\ref{condbeta}) is satisfied by $\alpha^2|\gvct{q}|^4/\beta^2$
or $\alpha^2|\gvct{q}|^4/(2\beta^2)$ of them, respectively. Thus
$\gvct{L}_P/\gvct{L}_{\gvct{q}}$ contains at most
$\alpha^4|\gvct{q}|^4/\beta^2$
vectors satisfying condition~(\ref{condbeta}) and hence $n'_P$ is a divisor
of $\alpha^4|\gvct{q}|^4/\beta^2$. Let $\gvct{L}^\beta_P$ denote the subgroup
of $\gvct{L}_P/\gvct{L}_{\gvct{q}}$ that is formed by the vectors satisfying
cond.~(\ref{condbeta}) and assume $\gvct{x}\in\gvct{L}^\beta_P$ in the
following. We can rewrite \Eq{condRx} as
\BAL
R(\gvct{p},\gvct{q})\gvct{x}=\frac{1}{|\gvct{p}\gvct{q}|}
\sum_{i=0}^{3}\gvct{u}_i\ip{\gvct{x}}{\bar{\gvct{p}}\gvct{u}_i\gvct{q}}=
\sum_{i=0}^{3}\gvct{u}_i
\frac{\ip{\gvct{x}}{\bar{\gvct{p}}\gvct{u}_i\gvct{q}}}{\alpha\beta\gamma},
\EAL
i.e. $|\gvct{p}\gvct{q}|$ must divide
$\ip{\gvct{x}}{\bar{\gvct{p}}\gvct{u}_i\gvct{q}}$
By assumption, $\beta$ divides
$\ip{\gvct{x}}{\bar{\gvct{p}}\gvct{u}_i\gvct{q}}$. On the other hand,
since $\den(R)=|\gvct{p}\gvct{q}|/2^\ell$, there exists an element $\gvct{x}$
of order $\alpha\gamma/2^\ell$ or higher. Thus at most
$2^\ell|\gvct{L}^\beta_P|/(\alpha\gamma)$ vectors
$\gvct{x}\in\gvct{L}_P/\gvct{L}_{\gvct{q}}$ satisfy condition~(\ref{condRx})
and hence $n'_P$ divides $2^{\ell+1}\alpha^3|\gvct{q}|^4/(\beta^2\gamma)=
2^{\ell+1}\alpha^3|\gvct{q}|^2$. From this we infer that $\Sigma_F(R)$
is a multiple of $\alpha|\gvct{q}|^2/2^{\ell+2}$ and hence a multiple of
$\Sigma(\gvct{q})$. Analogously we prove that $\Sigma_F(R)$
is a multiple of $\Sigma(\gvct{p})$. Thus
$\lcm(\Sigma(\gvct{p}),\Sigma(\gvct{q}))\leq \Sigma_F(R)$ and the claim
follows. \hfill$\Box$

From this result we can easily infer the coincidence index $\Sigma_P(R)$ for
the primitive lattice. Since $\gvct{L}_F$ is a sublattice of index $2$ of
$\gvct{L}_P$, $\Sigma_P(R)$ must divide $2\Sigma_F(R)$ and $\Sigma_F(R)$
must divide $2\Sigma_P(R)$.\cite{baa97} Since $\Sigma_F(R)$ is odd we have
$\Sigma_P(R)=\Sigma_F(R)$ or $\Sigma_P(R)=2\Sigma_F(R)$. Due to \Eq{dens}
the index $\Sigma_P(R)$ is odd if $\den(R)$ is odd and even if $\den(R)$ is
even. Hence we have proved
\begin{theo}
Let $\gvct{p},\gvct{q}$ be an admissible pair of primitive integer
quaternions and let $\Sigma(\gvct{q})=2^{-\ell}|\gvct{q}|^2$,
where $\ell=0,1,2$ is the maximal power such that $2^\ell$ divides
$|\gvct{q}|^2$. Then, for the primitive lattice,
the rotation $R(\gvct{p},\gvct{q})$ has coincidence
index
\BAL
\Sigma_P(\gvct{p},\gvct{q}):=\Sigma_P(R(\gvct{p},\gvct{q}))=
\lcm[\Sigma(\gvct{p}),\Sigma(\gvct{q}),\den(R(\gvct{p},\gvct{q}))].
\EAL
\end{theo}
This was first stated, without proof, in~\cite{baa97}.

\section{Equivalent CSLs}

Different coincidence rotations may generate the same CSL or rotated
copies of each other.
It is natural to group these rotations and CSLs in appropriately
chosen equivalence classes. The natural way is to call two coincidence
rotations \emph{equivalent}
if they are in the same double coset of the symmetry group of the
lattice~\cite{pzcsl1,grim74,grim76}. To be precise,
let $G_P$ and $G_F$ denote the symmetry groups
of the primitive and the face--centered hypercubic lattice. Then we call
two coincidence rotations $R$, $R'$ \emph{\mbox{P--e}quivalent}
(\emph{\mbox{F--e}quivalent}) if there
exist two rotations $Q, Q'\in G_P$ ($Q, Q'\in G_F$) such that $R=QR'Q'$.
Accordingly, we call two CSLs \mbox{P--e}quivalent
(\mbox{F--e}quivalent) if the
corresponding coincidence rotations are \mbox{P--e}quivalent
(\mbox{F--e}quivalent).
In particular, $R$ and $RQ$, $Q\in G_{P,F}$ give rise to the same CSL.

Hence two coincidence rotations are equivalent if they belong to the same
double coset $G_P R G_P$ or $G_F R G_F$. These double cosets can be
calculated if one knows the subgroups $H_i(R):=G_i \cap R G_i R^{-1}$,
$i=P,F$. In order to determine these groups we make use of the fact
that $SU(2)\times SU(2)$ is a double cover of the $4$--dimensional
rotation group $SO(4)$, which is reflected in the
parameterization~\Eq{parrot}. Although the corresponding double
cover of $G_F$ and $G_P$ is not a direct product but a subdirect product,
we can make use of this special property and reduce the $4$--dimensional
case to the $3$--dimensional one.

In order to do this we recall that the $3$--dimensional rotations can be
parameterized by quaternions as well~\cite{koecheng,hurw,val}. The group
$\CG$ of order $|\CG|=48$ generated by the quaternions $(\pm 1,0,0,0)$,
$\frac{1}{\sqrt{2}}(\pm 1,\pm 1,0,0)$, $\frac{1}{2}(\pm 1,\pm 1,\pm 1,\pm 1)$
and permutations thereof
is a double cover of the cubic symmetry group $O$ of order $|O|=24$.
Based on the notion of equivalence of $3$--dimensional coincidence rotations
we introduce the following equivalence notion for
quaternions: Two quaternions $\gvct{q}$ and $\gvct{q}'$ are \emph{equivalent} 
($\gvct{q}\sim\gvct{q}'$) if
there exist quaternions $\gvct{s},\gvct{s}'\in\CG$ such that
$\gvct{q}'=\gvct{s}\gvct{q}\gvct{s}'$. Their equivalence classes are
known~\cite{pzcsl1} and the different types
are summarized in Table~\ref{tab1}. Here
$\CH(\gvct{q}):=\CG\cap \gvct{q} \CG \gvct{q}^{-1}$. Furthermore the number
of inequivalent CSLs for a given $\Sigma$ is known~\cite{grim76,pzcsl1}.
These numbers are summarized in Table~\ref{tab3} for all special
quaternions $\gvct{q}$. The number of inequivalent CSLs for a general
$\gvct{q}$ can be obtained by considering the total number of
CSLs~\cite{pzcsl1}.

Let $\CG'\subseteq\CG$ be the group generated by the quaternions
$(\pm 1,0,0,0)$, $\frac{1}{2}(\pm 1,\pm 1,\pm 1,\pm 1)$
and permutations thereof. Now the group
$\CG_F=\CG'\otimes\CG'\cup
(\frac{1}{\sqrt{2}}(1,1,0,0),\frac{1}{\sqrt{2}}(1,1,0,0))\CG'\otimes\CG'$
is a double cover of $G_F$. We call two pairs of quaternions
$(\gvct{p},\gvct{q})$ and $(\gvct{p}',\gvct{q}')$ \emph{\mbox{F--e}quivalent}
if the
corresponding rotations $R(\gvct{p},\gvct{q})$ and $R(\gvct{p}',\gvct{q}')$
are \mbox{F--e}quivalent. If $(\gvct{p},\gvct{q})$ and
$(\gvct{p}',\gvct{q}')$ are
\mbox{F--e}quivalent then $\gvct{p}\sim\gvct{p}'$ and $\gvct{q}\sim\gvct{q}'$,
but the converse is not true in general. Let us analyze the converse
situation. Let $\gvct{p}\sim\gvct{p}'$ and $\gvct{q}\sim\gvct{q}'$, i.e.
there exist $\gvct{r},\gvct{r}',\gvct{s},\gvct{s}'\in\CG$ such that 
$\gvct{p}'=\gvct{r}\gvct{p}\gvct{r}'$ and
$\gvct{q}'=\gvct{s}\gvct{q}\gvct{s}'$. If both pairs $(\gvct{p},\gvct{q})$
and $(\gvct{p}',\gvct{q}')$ are admissible, then
$(\gvct{r}\gvct{r}',\gvct{s}\gvct{s}')$ must be admissible, too. If
$(\gvct{r},\gvct{s})$ is admissible, then so is $(\gvct{r}',\gvct{s}')$,
and $(\gvct{p},\gvct{q})$ and $(\gvct{p}',\gvct{q}')$ are \mbox{F--e}quivalent.
If $(\gvct{r},\gvct{s})$ is \emph{not} admissible,
then $(\gvct{p},\gvct{q})$ and
$(\gvct{p}',\gvct{q}')$ are \mbox{F--e}quivalent only if there exist
admissible pairs $(\gvct{r}_1,\gvct{s}_1)$ and $(\gvct{r}'_1,\gvct{s}'_1)$
such that $\gvct{r}\gvct{p}\gvct{r}'=\gvct{r}_1\gvct{p}\gvct{r}'_1$ and
$\gvct{s}\gvct{q}\gvct{s}'=\gvct{s}_1\gvct{q}\gvct{s}'_1$. This is possible if
and only if $\CH(\gvct{p})$ or $\CH(\gvct{q})$ contains one of the
quaternions $\frac{1}{\sqrt{2}}(\pm 1,\pm 1,0,0)$ or a permutation thereof.
The latter statement is equivalent to the statement that $\gvct{p}$ or
$\gvct{q}$ is equivalent to one of the following quaternions: $(1,0,0,0)$,
$(0,1,1,1)$, $(m,n,0,0)$ or $(m,n,n,0)$. We can summarize these
considerations as follows: If $(\gvct{p},\gvct{q})$ and
$(\gvct{p}',\gvct{q}')$ are
\mbox{F--e}quivalent then $\gvct{p}\sim\gvct{p}'$ and $\gvct{q}\sim\gvct{q}'$.
Conversely $\gvct{p}\sim\gvct{p}'$ and $\gvct{q}\sim\gvct{q}'$ implies that
$(\gvct{p},\gvct{q})$ and $(\gvct{p}',\gvct{q}')$ are \mbox{F--e}quivalent if
$\gvct{p}$ or $\gvct{q}$ is equivalent to one of the following quaternions:
$(1,0,0,0)$, $(0,1,1,1)$, $(m,n,0,0)$ or $(m,n,n,0)$.

Assume now that $\gvct{p}\sim (m,n,n,n)$ and $\gvct{q}\sim (m',n',n',n')$.
Then we may not conclude that the admissible pair $(\gvct{p},\gvct{q})$ is
\mbox{F--e}quivalent to $((m,n,n,n),(m',n',n',n'))$. However, we may conclude
that $(\gvct{p},\gvct{q})$ is \mbox{F--e}quivalent either to
$((m,n,n,n),(m',n',n',n'))$ or $((m,n,n,n),(m',-n',-n',-n'))$. Note that the
latter pairs are not \mbox{F--e}quivalent. Nevertheless, they are of
the same type.

Having this in mind we can use Table~\ref{tab1} to calculate all
types of possible \mbox{F--e}quivalence classes. Instead of calculating the
groups $H_F(R(\gvct{p},\gvct{q}))$ directly we compute their corresponding
double covers $\CH_F(\gvct{p},\gvct{q})$. It turns out that they are simply
given by
$\CH_F(\gvct{p},\gvct{q})=(\CH(\gvct{p})\otimes\CH(\gvct{q}))\cap\CG'$.
The results are listed in
Table~\ref{tab2}. In order to save space we have omitted some pairs. These
can be easily obtained by interchanging the role of $\gvct{p}$ and $\gvct{q}$
and adapting the corresponding subgroup $\CH_F$. In addition,
we have used the definition $\CH'_i:=\CH_i\cap\CG'$.

The fact that $\CG_F$ is a special subgroup of $\CG\times\CG$ enables us
to derive the number of different and inequivalent CSLs from the
$3$--dimensional case. First, we consider the total number of different CSLs
$f_F(\Sigma)$.
Recall that the total number of different CSLs $f(\Sigma)$ for a given
$\Sigma$ in the $3$--dimensional case is given by~\cite{pzcsl1,baa97}
\BAL
f(1)&=1\\
f(2)&=0\\
f(p^r)&=(p+1)p^{r-1}&& \mbox{if $p$ is an odd prime, $r\geq 1$}\\
f(mn)&=f(m)f(n) && \mbox{if $m,n$ are coprime.}\label{fsigma3}
\EAL
The multiplicativity of $f(\Sigma)$ is due to the uniqueness of the (left)
prime factorization of the integer quaternions~\cite{hurw}.
The same reasoning holds
true in four dimensions, too, so we only need to calculate $f_F(p^r)$.
To this end we note that there are precisely
$f(\Sigma(\gvct{p}))f(\Sigma(\gvct{q}))$ different CSLs
for given $\Sigma(\gvct{p})$ and $\Sigma(\gvct{q})$. Summing up all 
admissible combinations
of $(\Sigma(\gvct{p}),\Sigma(\gvct{q}))$ that give a fixed $p^r$ we
obtain~\cite{baa97}
\BAL\label{totalF}
f_F(p^r)
&=\frac{p+1}{p-1}p^{r-1}(p^{r+1}+p^{r-1}-2).
\EAL

In a similar way we can calculate the number of inequivalent CSLs of a certain
type, say $(\gvct{p},\gvct{q})\equiv ((0,1,1,1),(0,1,1,1))$ or
$((m,n,n,n),(m',n',n',n'))$. These results are summarized in
Tables~\ref{tab4a}~and~\ref{tab4b}.

Let us discuss some of them. Consider pairs of type $((1,0,0,0),(m,n,n,n))$
first. Then $\Sigma(\gvct{p},\gvct{q})=\Sigma_F$ implies
$\Sigma(\gvct{q})=\Sigma_F$. Hence the number of inequivalent CSLs is equal to
the number of inequivalent quaternions $\gvct{q}=(m,n,n,n)$, which can
be read off directly from Table~\ref{tab3}. Thus there are precisely $2^{k-1}$
inequivalent CSLs if $p=1\bmod 6$ for all the $k$ different prime factors of
$\Sigma_F$. Note that a prime factor $3$ cannot exist,
since $\Sigma_F$ must be a
square as $(\gvct{p},\gvct{q})$ must be an admissible pair.

Consider now pairs of type $((m,n,n,n),(m',n',n',n'))$. Such pairs can only
exist if
$\Sigma_F(\gvct{p},\gvct{q})=3^t\prod_{i}p_i^{2r_i}\prod_{j}q_j^{2s_j+1}$
where all prime factors $p_i,q_j=1\bmod 6$ and $t=0,1$. This implies that
$\Sigma(\gvct{p})=3^t\prod_{i}p_i^{2r'_i}\prod_{j}q_j^{2s'_j+1}$ and
$\Sigma(\gvct{q})=3^t\prod_{i}p_i^{2r''_i}\prod_{j}q_j^{2s''_j+1}$ with
$r_i=\max(r'_i,r''_i), s_j=\max(s'_j,s''_j)$. For a fixed combination
of $\{r'_i,s'_j\}$, there are $2^{k'-1}$ inequivalent quaternions $\gvct{p}$,
where $k'$ is the number of different prime factors $p_i,q_j\ne 3$ contained in
$\Sigma(\gvct{p})$. If we use the notation $\nu(a)=0,1$ for $a=0$, $a\geq 1$,
respectively, we can write
$2^{m'-1}=1/2 \prod_i 2^{\nu(2r'_i)}\prod_j 2^{\nu(2s'_j+1)}$. In order to
get the number of inequivalent admissible pairs $(\gvct{p},\gvct{q})$
we have to take the sum over all possible
combinations of $r'_i,r''_i,s'_j,s''_j$. Note that $r'_i$ runs through
$0,\ldots, r_i$ if $r''_i=r_i$ and vice versa. Hence the number of
inequivalent admissible pairs reads
\BAL
\lefteqn{\frac{1}{2}\sum_{(r'_i,r''_i,s'_j,s''_j)}
\prod_i 2^{\nu(2r'_i)+\nu(2r''_i)}
\prod_j 2^{\nu(2s'_j+1)+\nu(2s''_j+1)}}\\
&=
\frac{1}{2}\prod_i \left(\sum_{r'_i=0}^{r_i-1} 2^{\nu(2r'_i)+\nu(2r_i)}+
\sum_{r''_i=0}^{r_i-1} 2^{\nu(2r_i)+\nu(2r''_i)}+2^{\nu(2r_i)+\nu(2r_i)}\right)
\\
&\qquad \cdot
\prod_j \left(\sum_{s'_j=0}^{s_j-1}2^{\nu(2s'_j+1)+\nu(2s_j+1)}+
\sum_{s''_j=0}^{s_j-1}2^{\nu(2s_j+1)+\nu(2s''_j+1)}+
2^{\nu(2s_j+1)+\nu(2s''_j+1)}\right)\\
&=\frac{1}{2}\prod_i (8r_i)\prod_j [4(2s_j+1)]=
2\cdot 4^{k-1}\prod_i (2r_i)\prod_j (2s_j+1),
\EAL
where $k$ is the number of different prime factors $p_i,q_j\ne 3$. 
If $\Sigma_F(\gvct{p},\gvct{q})$ contains at least one
odd prime power $q_j^{2s_j+1}$, 
we have finished. Otherwise we have to take into account that the sum above
includes $2^{k-1}$ pairs of the form $((m,n,n,n),(1,0,0,0))$ or
$((m,n,n,n),(0,1,1,1))$. Hence a term $2^k$ must be subtracted from the sum
above. Thus there exist
\BAL
n_{F22}=2\cdot 4^{k-1}\prod_\ell t_\ell - \delta\, 2^{k} 
\EAL
inequivalent admissible pairs $((m,n,n,n),(m',n',n',n'))$ for a fixed
$\Sigma_F(\gvct{p},\gvct{q})=3^t \prod_{\ell}p_\ell^{t_\ell}$ with
$p_\ell=1\bmod 6$ and $t=0,1$. Here $\delta=1$ if all $t_\ell$ are even and
$\delta=0$ otherwise.

Next we consider the case $((m,n,n,n),(m',n',0,0))$, which is an example where
$\gvct{p}$ and $\gvct{q}$ are of different type. First observe that
$\Sigma(\gvct{q})$ may only contain prime factors $p=1 \bmod 4$, whereas
$\Sigma(\gvct{p})$ may only contain prime factors $p=1 \bmod 6$ and $p=3$,
for the latter only the powers $3^0$ and $3^1$ are allowed. Since the pair
must be admissible, the factor $p=3$ is ruled out and
the coincidence index takes the form
$\Sigma_F(\gvct{p},\gvct{q})=
\prod_i p_i^{2r_i}\prod_j p'_j{}^{2r'_j}\prod_\ell q_\ell^{s_\ell}$
where $p_i= 1\bmod 4, p_i \ne 1\bmod 6$,
$p'_j= 1\bmod 6, p'_j\ne 1\bmod 4$, $q_\ell=1\bmod 4, q_\ell=1\bmod 6$.
Hence $\Sigma(\gvct{q})=\prod_i p_i^{2r_i}\prod_\ell q_\ell^{s'_\ell}$,
$\Sigma(\gvct{p})=\prod_j p'_j{}^{2r'_j}\prod_\ell q_\ell^{s''_\ell}$ where
$s_\ell=\max(s'_\ell,s''_\ell)$. Again we have to sum over all possible
combinations $s'_\ell,s''_\ell$ and finally obtain the number $n_{F23}$
of \mbox{F--i}nequivalent admissible pairs
\BAL
n_{F23}=2^{k_1+k_2}4^{k_3-1}\prod_\ell s_\ell,
\EAL
if $k_1\geq 1$ and $k_2\geq 1$. Here $k_1, k_2, k_3$ are the number of
different prime factors $p_i$, $p'_j$, $q_\ell$. If $k_1=0,k_2\neq 0$
this expression
includes the pairs of type $((m,n,n,n),(1,0,0,0))$, so that a term
$2^{k_2+k_3-1}$ must be subtracted. Thus
\BAL
n_{F23}=2^{k_2}(4^{k_3-1}\prod_\ell s_\ell-2^{k_3-2}).
\EAL
A similar expression is obtained for $k_2=0$. Finally, if $k_1=k_2=0$, we
get
\BAL
n_{F23}= 4^{k_3-1}\prod_\ell s_\ell-2^{k_3-1}.
\EAL

At last, let us consider pairs where at least one quaternion is completely
general. As an example, we use $((m,n,n,n),\gvct{q})$. In this case, the
approach is slightly different from the previous cases, since we lack a nice
formula for the three--dimensional case. But we can proceed as follows: We
first calculate the number of different admissible pairs
$((m,n,n,n),\gvct{q})$, where $\gvct{q}$ is a general or a special quaternion.
We then subtract the number of all special combinations $((m,n,n,n),\gvct{q})$
and finally divide by the number of equivalent pairs.
We first note that $\Sigma((m,n,n,n),\gvct{q})$ must be of the form
$\Sigma=3^r\prod_{i}p_i^{s_i}\prod_jq_j^{2t_j}$, where $p_i=1\pmod 6$ and
$q_j\ne 1\pmod 6$ and $r\geq 0$ and at least one $s_i\geq 1$.
We have to sum over all pairs with
$\Sigma(m,n,n,n)=3^{r'}\prod_{i}p_i^{s'_i}$,
$\Sigma(\gvct{q})=3^r\prod_{i}p_i^{s''_i}\prod_jq_j^{2t_j}$ such that
$r'\leq 1, r'=r\pmod2, s_i=\max(s'_i,s''_i)$.
For fixed $\Sigma(m,n,n,n)$ and $\Sigma(\gvct{q})$ we have the following
situation: There are $2^{k-1}=1/2\prod_{i} 2^{1-\delta_{0,s'_i}}$ inequivalent
quaternions of type $(m,n,n,n)$ (if at least one $s'_i>0$, $k$ is the number
of different prime factors $>3$) and there are
$48\cdot(4\cdot 3^{r-1})^{1-\delta_{0,r}}\prod_{i}(p_i+1)p_i^{s''_i-1}\prod_j(q_j+1)q_j^{2t_j-1}$
different (in general \emph{not} inequivalent) quaternions $\gvct{q}$.
Note that the product
ranges only over those $i$ for which $s''_i>0$. If we use Gauss' symbol
$[x]$ in order to denote the largest integer $n\leq x$ we may rewrite
this as
$48\cdot[4\cdot 3^{r-1}]\prod_{i}[(p_i+1)p_i^{s''_i-1}]\prod_j(q_j+1)q_j^{2t_j-1}$
and take the product over all $i$. Hence for fixed $\Sigma(m,n,n,n)$ and
$\Sigma(\gvct{q})$ we have
\BAL
1/2\cdot 8\cdot 48 \cdot 1/2 \prod_{i} 2^{1-\delta_{0,s'_i}}\cdot
 48\cdot[4\cdot 3^{r-1}]\prod_{i}[(p_i+1)p_i^{s''_i-1}]\prod_j(q_j+1)q_j^{2t_j-1}
\EAL
different (in general \emph{not} inequivalent admissible pairs). Note that we
have added a factor $1/2$ taking into account that only half of the pairs
are admissible. Summing over all possible combinations of $\Sigma(m,n,n,n)$
and $\Sigma(\gvct{q})$ we get
\BAL
1152m_{F2}&=4\cdot1152\cdot[4\cdot 3^{r-1}]\prod_{i}
\left(\sum_{\ell_i=1}^{[s_i/2]}2^{1-\delta_{0,s_i-2\ell_i}}(p_i+1)p_i^{s_i-1}+
\sum_{\ell_i=0}^{[s_i/2]} 2[(p_i+1)p_i^{s_i-2\ell_i-1}]\right)
\prod_j(q_j+1)q_j^{2t_j-1}\\
&=
4\cdot1152\cdot[4\cdot 3^{r-1}]\prod_{i}
\left((s_i+1)(p_i+1)p_i^{s_i-1}+2\frac{p_i^{s_i-1}-1}{p_i-1}\right)
\prod_j(q_j+1)q_j^{2t_j-1}
\EAL
different admissible pairs if there is at least one $s_i$ is odd. Otherwise
we must exclude the term with $\ell_i=s_i/2$ for all $i$ in the first sum,
i.e.
\BAL
m_{F2}&=4\cdot[4\cdot 3^{r-1}]\left(\prod_{i}
\left((s_i+1)(p_i+1)p_i^{s_i-1}+2\frac{p_i^{s_i-1}-1}{p_i-1}\right)-
\delta_{F2}\prod_{i}(p_i+1)p_i^{s_i-1}\right)
\prod_j(q_j+1)q_j^{2t_j-1},
\EAL
where $\delta_{F2}=0,1$ according to whether there exists an odd $s_i$ or not.
From this expression we subtract all admissible pairs with special $\gvct{q}$,
divide by the number of equivalent pairs and obtain the following expression
for the number of inequivalent admissible quaternions of type
$((m,n,n,n),\gvct{q})$:
\BAL
n_{F25}=\frac{1}{96}\left(m_{F2}-\sum_{i=0}^4 g_{F2i}n_{F2i}\right).
\EAL
Similar expression are obtained for $n_{F35}$ and $n_{F45}$. And finally
we can compute $n_{F55}$ by recalling the total number of different
quaternions $f_F$ given in \Eq{totalF}:
\BAL
f_F=\sum_{i,j=0}^5 g_{Fij}n_{Fij}.
\EAL

Finally let us have a short look on the primitive hypercubic lattice.
Similar results can be proved for this case. The best way to obtain them
is to derive them directly from the previous results. We just have to keep in
mind that the symmetry group $G_P$ is a subgroup of index 3 of $G_F$. In
particular, the coset decomposition for the corresponding groups of
quaternions reads
\BAL
\CG_F=\CG_P\cup\left((1,0,0,0),\textstyle{\frac{1}{2}}(1,1,1,1)\right)\CG_P
\cup \left((1,0,0,0),\textstyle{\frac{1}{2}}(1,-1,-1,-1)\right)\CG_P.
\EAL
If we apply this decomposition to the double cosets
$\CG_F(\gvct{p},\gvct{q})\CG_F$, we get the double cosets of $\CG_P$,
which are just the \mbox{P--e}quivalence classes of admissible pairs,
see~Tab.~\ref{tabefb}. The corresponding groups $\CH_P(\gvct{p},\gvct{q})$
can now be inferred from the corresponding groups $\CH_F(\gvct{p},\gvct{q})$.
In particular, we have
$\CH_P(\gvct{p},\gvct{q})\subseteq\CH_F(\gvct{p},\gvct{q})\cap\CG_P$, which
simplifies the determination of $\CH_P(\gvct{p},\gvct{q})$ considerably.
The results are shown in Tab.~\ref{tabebpartone}. 
Combining these results with the numbers $n_{Fij}$ of \mbox{F--i}nequivalent 
admissible pairs, we get the number of \mbox{P--i}nequivalent admissible pairs,
which are listed in Tab.~\ref{tabip}.

\section{Conclusions and Outlook}

We have calculated the coincidence index $\Sigma$ for both kinds of 
four--dimensional hypercubic
lattices. Moreover, we have determined all CSLs and their equivalence
classes as
well as the total number of different and inequivalent CSLs for fixed
$\Sigma$. Here, equivalence always means equivalence up to proper rotations.
But of course there exist
reflections that leave the hypercubic lattices invariant and one can be
interested in extending the notion of equivalence to the full symmetry group.
We briefly sketch
how one can include the improper rotations.
First note that the special reflection
$m:\gvct{q}\to (q_0,-q_1,-q_2,-q_3)$  just corresponds to quaternion
conjugation. Now any symmetry operation is a product of this reflection and a
rotation, and it is sufficient to consider this reflection in detail.
Since $mR(\gvct{p},\gvct{q})=R(\gvct{q},\gvct{p})m$, it follows that
the admissible pairs $(\gvct{p},\gvct{q})$ and $(\gvct{q},\gvct{p})$ are
equivalent. Thus we have two situations: If $(\gvct{p},\gvct{q})$ and
$(\gvct{q},\gvct{p})$ are not equivalent under proper rotations, than their
equivalence classes merge to form a single equivalence class. If
$(\gvct{p},\gvct{q})$ and $(\gvct{q},\gvct{p})$ are already equivalent under
proper rotations, than the equivalence class stays the same and the
corresponding symmetry group $H(\gvct{p},\gvct{q})$ contains a symmetry
operation which is a conjugate of $m$. Thus we know all
equivalence classes and their symmetry groups $H(\gvct{p},\gvct{q})$. It is
then straightforward to calculate the number of inequivalent CSLs.

\section*{Acknowledgements}

The author is very grateful to Michael Baake for interesting discussions
on the present subject and to the Faculty of Mathematics, University
Bielefeld, for its hospitality. Financial support by the Austrian
Academy of Sciences (APART-program) is gratefully acknowledged.

\begin{table}[htb]
\begin{center}
\begin{tabular}{|c|c|c|c|}
\hline
$\gvct{q}$ & $\CH(\gvct{q})$ & $|\CH(\gvct{q})|$ & $|\CG\gvct{q}\CG|$ \\
\hline
$(1,0,0,0)$ & $\CG$ & 48 & 48 \\
\hline
$(0,1,1,1)\sim(3,1,1,1)$ & $\CH_1=[(-1,0,0,0),(1,1,1,1),(0,1,-1,0)]$ & $12$ &
$4\cdot 48=192$ \\
\hline
$(m,n,n,n)$ & $\CH_2=[(-1,0,0,0),(1,1,1,1)]$ & $6$ & $8\cdot 48=384$ \\
\hline
$(m,n,0,0)$ & $\CH_3=[(-1,0,0,0),(1,1,0,0)]$ & $8$ & $6\cdot 48=288$ \\
\hline
$(m,n,n,0)$ & $\CH_4=[(-1,0,0,0),(0,1,1,0)]$ & $4$ & $12\cdot 48=576$ \\
\hline
otherwise & $\CH_5=[(-1,0,0,0)]$ & $2$ & $24\cdot 48=1152$ \\
\hline
\end{tabular}
\end{center}
\caption{\label{tab1} Equivalence classes of quaternions: Any primitive 
quaternion is equivalent to one of the quaternions in the first column.
The second column gives a set of generators of $\CH(\gvct{q})$. The third
column gives the order of $|\CH(\gvct{q})|$ and the forth column states
the number of equivalent $\gvct{q}$, which is 48 times the number of
equivalent $3$--dimensional CSLs.}
\end{table}

\begin{table}[htb]
\begin{center}
\begin{tabular}{|c|c|c|c|}
\hline
$\gvct{q}$ & inequiv. CSLs & condition \\
\hline
$(1,0,0,0)$ & $1$ & $\Sigma=1$ \\
\hline
$(0,1,1,1)$ & $1$ & $\Sigma=3$ \\
\hline
$(m,n,n,n)$ & $2^{k-1}$ &
\begin{minipage}{11cm}
$p=1\bmod 6$ for all prime factors $p\ne 3$ of
$\Sigma>3$, the factor $p=3$ occurs at most once and $k$ is the number of
different prime factors $p=1\bmod 6$ of $\Sigma$
\end{minipage} \\
\hline
$(m,n,0,0)$ & $2^{k-1}$ &
\begin{minipage}{11cm}
$p=1\bmod 4$ for all prime factors $p$ of
$\Sigma$ and $k$ is the number of different prime factors of $\Sigma$.
\end{minipage} \\
\hline
$(m,n,n,0)$ & $2^{k-1}$ &
\begin{minipage}{11cm}
 $p=1$ or $3\bmod 8$ for all prime factors $p$ of
$\Sigma$, where $k$ is the number of different prime factors of $\Sigma>3$.
\end{minipage} \\
\hline
\end{tabular}
\end{center}
\caption{\label{tab3} Number of inequivalent cubic CSLs/coincidence rotations
for a fixed value $\Sigma$. The last column gives the condition under which
these values hold. If this condition is not satisfied, the corresponding
number of inequivalent CSLs is 0 for the particular type of $\gvct{q}$.}
\end{table}

\begin{table}[htb]
\begin{center}
\begin{tabular}{|c|c|c|c|c|}
\hline
$\gvct{p}$ & $\gvct{q}$ & $\CH_F(\gvct{p},\gvct{q})$ &
$|\CH_F(\gvct{p},\gvct{q})|$ & $|G_F R(\gvct{p},\gvct{q}) G_F|$ \\
\hline
$(1,0,0,0)$ & $(1,0,0,0)$ & $\CG_F$ & $1152$ & $576g_{F00}=576$ \\
\hline
$(1,0,0,0)$ & $(m,n,n,n)$ & $\CG'\otimes\CH'_2$ &
$144$ & $576g_{F02}=8\cdot 576$ \\
\hline
$(1,0,0,0)$ & $(m,n,0,0)$ & $\CG'\otimes\CH'_3\cup
(\frac{1}{\sqrt{2}}(1,1,0,0),\frac{1}{\sqrt{2}}(1,1,0,0))\CG'\otimes\CH'_3$ &
$192$ & $576g_{F03}=6\cdot 576$\\
\hline
$(1,0,0,0)$ & $(m,n,n,0)$ & $\CG'\otimes\CH'_4\cup
(\frac{1}{\sqrt{2}}(1,1,0,0),\frac{1}{\sqrt{2}}(0,1,1,0))\CG'\otimes\CH'_4$ &
$96$ & $576g_{F04}=12 \cdot 576$\\
\hline
$(1,0,0,0)$ & general & $\CG'\otimes\CH'_5$ & $48$ & $576g_{F05}=24\cdot 576$\\
\hline
$(0,1,1,1)$ & $(0,1,1,1)$ & $\CH'_1\otimes\CH'_1\cup
(\frac{1}{\sqrt{2}}(0,1,-1,0),\frac{1}{\sqrt{2}}(0,1,-1,0))\CH'_1\otimes\CH'_1$
& $72$ & $576g_{F11}=16\cdot 576$\\
\hline
$(0,1,1,1)$ & $(m,n,n,n)$ & $\CH'_1\otimes\CH'_2$ & $36$ & $576g_{F12}=32\cdot 576$\\
\hline
$(0,1,1,1)$ & $(m,n,0,0)$ & $\CH'_1\otimes\CH'_3\cup
(\frac{1}{\sqrt{2}}(0,1,-1,0),\frac{1}{\sqrt{2}}(1,1,0,0))\CH'_1\otimes\CH'_3$
& $48$ & $576g_{F13}=24\cdot 576$\\
\hline
$(0,1,1,1)$ & $(m,n,n,0)$ & $\CH'_1\otimes\CH'_4\cup
(\frac{1}{\sqrt{2}}(0,1,-1,0),\frac{1}{\sqrt{2}}(0,1,1,0))\CH'_1\otimes\CH'_4$
& $24$ & $576g_{F14}=48\cdot 576$\\
\hline
$(0,1,1,1)$ & general & $\CH'_1\otimes\CH'_5$ & $12$ & $576g_{F15}=96\cdot 576$\\
\hline
$(m,n,n,n)$ & $(m',n',n',n')$ & $\CH'_2\otimes\CH'_2$ & $36$ & $576g_{F22}=32\cdot 576$\\
\hline
$(m,n,n,n)$ & $(m',n',0,0)$ & $\CH'_2\otimes\CH'_3$ & $24$ & $576g_{F23}=48\cdot 576$\\
\hline
$(m,n,n,n)$ & $(m',n',n',0)$ & $\CH'_2\otimes\CH'_4$ & $12$ & $576g_{F24}=96\cdot 576$\\
\hline
$(m,n,n,n)$ & general & $\CH'_2\otimes\CH'_5$ & $12$ & $576g_{F25}=96\cdot 576$\\
\hline
$(m,n,0,0)$ & $(m',n',0,0)$ & $\CH'_3\otimes\CH'_3\cup
(\frac{1}{\sqrt{2}}(1,1,0,0),\frac{1}{\sqrt{2}}(1,1,0,0))\CH'_3\otimes\CH'_3$
& $32$ & $576g_{F33}=36\cdot 576$\\
\hline
$(m,n,0,0)$ & $(m',n',n',0)$ & $\CH'_3\otimes\CH'_4\cup
(\frac{1}{\sqrt{2}}(1,1,0,0),\frac{1}{\sqrt{2}}(0,1,1,0))\CH'_3\otimes\CH'_4$
& $16$ & $576g_{F34}=72\cdot 576$\\
\hline
$(m,n,0,0)$ & general & $\CH'_3\otimes\CH'_5$ & $8$ & $576g_{F35}=144\cdot 576$\\
\hline
$(m,n,n,0)$ & $(m',n',n',0)$ & $\CH'_4\otimes\CH'_4\cup
(\frac{1}{\sqrt{2}}(0,1,1,0),\frac{1}{\sqrt{2}}(0,1,1,0))\CH'_4\otimes\CH'_4$
& $8$ & $576g_{F44}=144\cdot 576$\\
\hline
$(m,n,n,0)$ & general & $\CH'_4\otimes\CH'_5$ & $4$ & $576g_{F45}=288\cdot 576$\\
\hline
general & general & $\CH'_5\otimes\CH'_5$ & $4$ & $576g_{F55}=288\cdot 576$\\
\hline
\end{tabular}
\caption{\label{tab2} F--Equivalence classes of admissible pairs. For each
admissible pair the corresponding group $\CH_F(\gvct{p},\gvct{q})$ and its
order is listed. The last column gives the number of equivalent coincidence
rotations $R(\gvct{p},\gvct{q})$. By dividing these numbers by 576 we obtain
the number of equivalent CSLs.
In order to save space we have omitted some pairs. These
can be easily obtained by interchanging the role of $\gvct{p}$ and $\gvct{q}$
and by adapting the subgroup $\CH_F$ correspondingly.}
\end{center}
\end{table}

\begin{table}[tb]
\begin{center}
\begin{tabular}{|c|c|c|c|}
\hline
$\gvct{p}$ & $\gvct{q}$ & inequivalent CSLs & condition\\
\hline
$(1,0,0,0)$ & $(1,0,0,0)$ & $n_{F00}=1$ & $\Sigma_F=1$\\
\hline
$(1,0,0,0)$ & $(m,n,n,n)$ & $n_{F02}=2^{k-1}$ & 
\begin{minipage}{6.1cm}
$k$ is the number of different prime
factors, all prime factors $p=1\bmod 6$, $\Sigma_F$ is a square
\end{minipage}\\
\hline
$(1,0,0,0)$ & $(m,n,0,0)$ & $n_{F03}=2^{k-1}$ &
\begin{minipage}{6.1cm}
$k$ is the number of different prime
factors, all prime factors $p=1\bmod 4$, $\Sigma_F$ is a square
\end{minipage}\\
\hline
$(1,0,0,0)$ & $(m,n,n,0)$ & $n_{F04}=2^{k-1}$ &
\begin{minipage}{6.1cm}
$k$ is the number of different prime
factors, all prime factors $p=1$ or $3\bmod 8$, $\Sigma_F$ is a square
\end{minipage}\\
\hline
$(1,0,0,0)$ & general & $n_{F05}$ & $\Sigma_F$ is a square\\
\hline
$(0,1,1,1)$ & $(0,1,1,1)$ & $n_{F11}=1$ & $\Sigma_F=3$\\
\hline
$(0,1,1,1)$ & $(m,n,n,n)$ & $n_{F12}=2^{k-1}$ &
\begin{minipage}{6.1cm}
$\Sigma_F=3a^2$, $k$ is the number of different prime
factors of $a$, all prime factors $p=1\bmod 6$
\end{minipage}\\
\hline
$(0,1,1,1)$ & $(m,n,n,0)$ & $n_{F14}=2^{k-1}$ &
\begin{minipage}{6.1cm}
$\Sigma_F=3a^2$, $k$ is the number of different prime
factors of $\Sigma_F$, all prime factors $p=1$ or $3\bmod 8$
\end{minipage}\\
\hline
$(0,1,1,1)$ & general & $n_{F15}$ & $\Sigma_F=3a^2$\\
\hline
$(m,n,n,n)$ & $(m',n',n',n')$ &
$n_{F22}=2\cdot 4^{k-1}\prod_\ell t_\ell - \delta\, 2^{k}$
&\begin{minipage}{6.1cm}
$\Sigma_F=3^{r}a$, $r=0,1$, $k$ is the number of different prime
factors of $a$, which is not divisible by $3$, all prime factors $p=1\bmod 6$,
$\delta=1$ if $a$ is a square and $\delta=0$ otherwise
\end{minipage}\\
\hline
$(m,n,n,n)$ & $(m',n',0,0)$ & 
 \begin{minipage}{6cm}
\begin{multline*}
n_{F23}=2^{k_1+k_2}4^{k_3-1}\prod_\ell t_\ell\\
-\delta_1 2^{k_2+k_3-2}-\delta_2 2^{k_1+k_3-2}
\end{multline*}
\end{minipage}
 & \begin{minipage}{6.1cm}$\Sigma_F=
\prod_i p_i^{2r_i}\prod_j p'_j{}^{2r'_j}\prod_\ell q_\ell^{t_\ell}$,
$p_i=1\pmod 4\ne 1\pmod 6$, $p'_j=1\pmod 6\ne 1\pmod 4$,
$q_\ell=1\pmod 4=1\pmod 6$, $k_1,k_2,k_3$ denote the number of different prime
factors of type $p_i$, $p'_j$ and $q_\ell$, respectively. $\delta_1=0$ unless
all $t_\ell$ are even and $k_1=0$, where $\delta_1=1$. An analogous definition
applies for $\delta_2$ with $k_1=0$ replaced by $k_2=0$.
\end{minipage}\\
\hline
$(m,n,n,n)$ & general & $n_{F25}$ & \\
\hline
\end{tabular}
\caption{\label{tab4a} Number of \mbox{F--i}nequivalent CSLs (Part 1).
The last column gives the condition under which
these values hold. If this condition is not satisfied, the corresponding
number of inequivalent CSLs is 0 for the particular type of
$(\gvct{p},\gvct{q})$. In order to save space we have omitted some pairs. These
can be easily obtained by interchanging $\gvct{p}$ and $\gvct{q}$ and reading
of the corresponding value $n_{Fij}=n_{Fji}$.}
\end{center}
\end{table}
\begin{table}[tb]
\begin{center}
\begin{tabular}{|c|c|c|c|}
\hline
$\gvct{p}$ & $\gvct{q}$ & inequivalent CSLs & condition\\
\hline
$(m,n,n,n)$ & $(m',n',n',0)$ &
\begin{minipage}{6cm}
\begin{multline*}
n_{F24}=2^{k_1+k_2}4^{k_3-1}\prod_\ell t_\ell\\
-1/2(n_{F02}+n_{F12}+n_{F04}+n_{F14})
\end{multline*}
\end{minipage} &
\begin{minipage}{6.1cm}
$\Sigma_F=
3^s\prod_i p_i^{2r_i}\prod_j p'_j{}^{2r'_j}\prod_\ell q_\ell^{t_\ell}$,
$p_i=1\pmod 6\ne 1\mbox{ or } 3\pmod 8$,
$p'_j=1\mbox{ or } 3\pmod 8\ne 1\pmod 6$,
$q_\ell=1\pmod 6=1\mbox{ or } 3\pmod 8$,
there must be at least one prime factor $=1\pmod 6$ and one $=1$ or $3\pmod 8$.
$k_1,k_2$ denote the number of different prime
factors of type $p_i$ and $p'_j$, respectively. $k_3$ is the number of prime
factors of type $q_\ell$ if $s=0$ and the number of prime
factors of type $q_\ell$ plus $1$  if $s>1$. 
\end{minipage}\\
\hline
$(m,n,0,0)$ & $(m',n',0,0)$ &
$n_{F33}=4^{k-1}\prod_\ell t_\ell - \delta\, 2^{k-1}$& 
\begin{minipage}{6.1cm}
$k$ is the number of different prime factors of $\Sigma_F$,
all prime factors $p=1\bmod 4$, $\delta=1$ if $\Sigma_F$ is a square and
$\delta=0$ otherwise
\end{minipage}\\
\hline
$(m,n,0,0)$ & $(m',n',n',0)$ &
\begin{minipage}{6cm}
\begin{multline*}
n_{F24}=2^{k_1+k_2}4^{k_3-1}\prod_\ell t_\ell\\
-1/2(n_{F02}+n_{F12}+n_{F04}+n_{F14}) 
\end{multline*}
\end{minipage}
& \begin{minipage}{6.1cm}
$\Sigma_F=
\prod_i p_i^{2r_i}\prod_j p'_j{}^{2r'_j}\prod_\ell q_\ell^{t_\ell}$,
$p_i=1\pmod 4\ne 1\mbox{ or } 3\pmod 8$,
$p'_j=3\pmod 8$,
$q_\ell=1\pmod 8$,
there must be at least one prime factor $=1\pmod 4$ and one $=1$ or $3\pmod 8$.
$k_1,k_2,k_3$ denote the number of different prime
factors of type $p_i$, $p'_j$, and $q_\ell$, respectively.
\end{minipage}\\
\hline
$(m,n,0,0)$ & general & $n_{F35}$ & \\
\hline
$(m,n,n,0)$ & $(m',n',n',0)$ &
$n_{F44}=4^{k-1}\prod_\ell t_\ell - \delta\, 2^{k-1}$&
\begin{minipage}{6.1cm}
$k$ is the number of different prime factors of $\Sigma_F$,
all prime factors $p=1,3\bmod 8$, $\delta=1$ if $\Sigma_F=a^2,3a^2$ and
$\delta=0$ otherwise
\end{minipage}\\
\hline
$(m,n,n,0)$ & general & $n_{F45}$ & \\
\hline
general & general & $n_{F55}$ & \\
\hline
\end{tabular}
\caption{\label{tab4b} Number of \mbox{F--i}nequivalent CSLs (Part 2)}
\end{center}
\end{table}

\begin{table}[htb]
\begin{center}
\begin{tabular}{|c|c|}
\hline
$\gvct{g}=(\gvct{p},\gvct{q})$&double coset decomposition of
$\CG_F\gvct{g}\CG_F$ \\
\hline
$((1,0,0,0),(1,0,0,0))$&$\CG_P\cup\gvct{s_1}\CG_P$\\
\hline
$((1,0,0,0),(m,n,n,n))$ & $\CG_P\gvct{g}\CG_P\cup
\CG_P\gvct{g}\gvct{s_1}\CG_P\cup\CG_P\gvct{g}\gvct{s_2}\CG_P$ \\
\hline
$((1,0,0,0),(m,n,0,0))$ & $\CG_P\gvct{g}\CG_P\cup\CG_P\gvct{g}\gvct{s_1}\CG_P$
 \\
\hline
$((1,0,0,0),(m,n,n,0))$ & $\CG_P\gvct{g}\CG_P\cup\CG_P\gvct{g}\gvct{s_1}\CG_P$ \\
\hline
$((1,0,0,0),(m,n,p,q))$ & $\CG_P\gvct{g}\CG_P\cup
\CG_P\gvct{g}\gvct{s_1}\CG_P\cup\CG_P\gvct{g}\gvct{s_2}\CG_P$ \\
\hline
$((0,1,1,1),(0,1,1,1))$ & $\CG_P\gvct{g}\CG_P\cup\CG_P\gvct{g}\gvct{s_1}\CG_P$ \\
\hline
$((0,1,1,1),(m,n,n,n))$ & $\CG_P\gvct{g}\CG_P\cup
\CG_P\gvct{g}\gvct{s_1}\CG_P\cup\CG_P\gvct{g}\gvct{s_2}\CG_P$ \\
\hline
$((0,1,1,1),(m,n,0,0))$ &$\CG_P\gvct{g}\CG_P\cup\CG_P\gvct{g}\gvct{s_1}\CG_P$  \\
\hline
$((0,1,1,1),(m,n,n,0))$ & $\CG_P\gvct{g}\CG_P\cup\CG_P\gvct{g}\gvct{s_1}\CG_P$ \\
\hline
$((0,1,1,1) ,(m,n,p,q))$ & $\CG_P\gvct{g}\CG_P\cup
\CG_P\gvct{g}\gvct{s_1}\CG_P\cup\CG_P\gvct{g}\gvct{s_2}\CG_P$  \\
\hline
$((m,n,n,n),(m',n',n',n'))$ & $\CG_P\gvct{g}\CG_P\cup
\CG_P\gvct{g}\gvct{s_1}\CG_P\cup\CG_P\gvct{g}\gvct{s_2}\CG_P$  \\
\hline
$((m,n,n,n),(m',n',0,0))$ & $\CG_P\gvct{g}\CG_P\cup
\CG_P\gvct{g}\gvct{s_1}\CG_P\cup\CG_P\gvct{g}\gvct{s_2}\CG_P$  \\
\hline
$((m,n,n,n),(m',n',n',0))$ & $\CG_P\gvct{g}\CG_P\cup
\CG_P\gvct{g}\gvct{s_1}\CG_P\cup\CG_P\gvct{g}\gvct{s_2}\CG_P$  \\
\hline
$((m,n,n,n),(m',n',p',q'))$ & $\CG_P\gvct{g}\CG_P\cup
\CG_P\gvct{g}\gvct{s_1}\CG_P\cup\CG_P\gvct{g}\gvct{s_2}\CG_P$  \\
\hline
$((m,n,0,0),(m',n',0,0))$ &
$\CG_P\gvct{g}\CG_P\cup\CG_P\gvct{g}\gvct{s_1}\CG_P\cup\CG_P\gvct{s_1}\gvct{g}\CG_P\cup\CG_P\gvct{s_1}\gvct{g}\gvct{s_1}\CG_P\cup\CG_P\gvct{s_1}\gvct{g}\gvct{s_2}\CG_P$  \\
\hline
$((m,n,0,0),(m',n',n',0))$ & $\CG_P\gvct{g}\CG_P\cup\CG_P\gvct{g}\gvct{s_1}\CG_P\cup\CG_P\gvct{g}\gvct{s_2}\CG_P\cup\CG_P\gvct{s_1}\gvct{g}\CG_P\cup\CG_P\gvct{s_1}\gvct{g}\gvct{s_2}\CG_P$ \\
\hline
$((m,n,0,0),(m',n',p',q'))$ & \begin{minipage}{9.5cm}$\CG_P\gvct{g}\CG_P\cup\CG_P\gvct{g}\gvct{s_1}\CG_P\cup\CG_P\gvct{g}\gvct{s_2}\CG_P\cup\CG_P\gvct{s_1}\gvct{g}\CG_P\cup\CG_P\gvct{s_2}\gvct{g}\CG_P\cup\CG_P\gvct{s_1}\gvct{g}\gvct{s_1}\CG_P\cup\CG_P\gvct{s_1}\gvct{g}\gvct{s_2}\CG_P\cup\CG_P\gvct{s_2}\gvct{g}\gvct{s_1}\CG_P\cup\CG_P\gvct{s_2}\gvct{g}\gvct{s_2}\CG_P$\end{minipage} \\
\hline
$((m,n,n,0),(m',n',n',0))$ & $\CG_P\gvct{g}\CG_P\cup\CG_P\gvct{g}\gvct{s_1}\CG_P\cup\CG_P\gvct{s_1}\gvct{g}\CG_P\cup\CG_P\gvct{s_1}\gvct{g}\gvct{s_1}\CG_P\cup\CG_P\gvct{s_1}\gvct{g}\gvct{s_2}\CG_P$ \\
\hline
$((m,n,n,0),(m',n',p',q'))$  &  \begin{minipage}{9.5cm}$\CG_P\gvct{g}\CG_P\cup\CG_P\gvct{g}\gvct{s_1}\CG_P\cup\CG_P\gvct{g}\gvct{s_2}\CG_P\cup\CG_P\gvct{s_1}\gvct{g}\CG_P\cup\CG_P\gvct{s_2}\gvct{g}\CG_P\cup\CG_P\gvct{s_1}\gvct{g}\gvct{s_1}\CG_P\cup\CG_P\gvct{s_1}\gvct{g}\gvct{s_2}\CG_P\cup\CG_P\gvct{s_2}\gvct{g}\gvct{s_1}\CG_P\cup\CG_P\gvct{s_2}\gvct{g}\gvct{s_2}\CG_P$\end{minipage}\\
\hline
$((m,n,p,q) ,(m',n',p',q'))$  &  \begin{minipage}{9.5cm}$\CG_P\gvct{g}\CG_P\cup\CG_P\gvct{g}\gvct{s_1}\CG_P\cup\CG_P\gvct{g}\gvct{s_2}\CG_P\cup\CG_P\gvct{s_1}\gvct{g}\CG_P\cup\CG_P\gvct{s_2}\gvct{g}\CG_P\cup\CG_P\gvct{s_1}\gvct{g}\gvct{s_1}\CG_P\cup\CG_P\gvct{s_1}\gvct{g}\gvct{s_2}\CG_P\cup\CG_P\gvct{s_2}\gvct{g}\gvct{s_1}\CG_P\cup\CG_P\gvct{s_2}\gvct{g}\gvct{s_2}\CG_P$\end{minipage}\\
\hline
\end{tabular}
\caption{\label{tabefb} Splitting of {F--e}quivalence classes into
{P--e}quivalence
  classes. The last column gives the decomposition of
the {F--e}quivalence class $\CG_F(\gvct{p},\gvct{q})\CG_F$ into double cosets
of $\CG_P$. Here, $\gvct{s_1}=(\gvct{u_0},\frac{1}{2}(1,1,1,1)),
 \gvct{s_2}=(\gvct{u_0},\frac{1}{2}(-1,1,1,1))$.
In order to save space we have omitted some pairs. These
can be easily obtained by interchanging $\gvct{p}$ and
$\gvct{q}$ and adapting the decomposition correspondingly, i.e. 
we have to interchange the corresponding quaternions of the pairs
$\gvct{g}\gvct{s_1}$, $\gvct{s_1}\gvct{g}$, \ldots as well.}

\end{center}
\end{table}

\begin{table}[htb]
\begin{center}
\begin{sideways}
\begin{tabular}{|c|c|c|c|c|}
\hline
$\gvct{p}$ & $\gvct{q}$ & non--trivial generators of
$\CH_P(\gvct{p},\gvct{q})$ &
$|\CH_P(\gvct{p},\gvct{q})|$ & $|G_P R(\gvct{p},\gvct{q}) G_P|$ \\
\hline
$(1,0,0,0)$ & $(1,0,0,0)$ & $\CG_P$ & $384$ & $192$ \\
$(1,0,0,0)$ & $(1,1,1,1)$ & $(\gvct{u}_1,\gvct{u}_0),(\gvct{u}_2,\gvct{u}_0),
(\gvct{u}_0,\gvct{u}_1),(\gvct{u}_0,\gvct{u}_2),(\frac{1}{2}(1,1,1,1),\frac{1}{2}(1,1,1,1))$ & $192$ & $2\cdot192$ \\
\hline
$(1,0,0,0)$ & $(m,n,n,n)$ & $(\gvct{u}_1,\gvct{u}_0),(\gvct{u}_2,\gvct{u}_0),
(\frac{1}{2}(1,1,1,1),\frac{1}{2}(1,1,1,1))$ &
$48$ & $8\cdot 192$ \\
\hline
$(1,0,0,0)$ & $(m,n,0,0)$ & $(\gvct{u}_1,\gvct{u}_0),(\gvct{u}_2,\gvct{u}_0),
(\frac{1}{\sqrt{2}}(1,1,0,0),\frac{1}{\sqrt{2}}(1,1,0,0))$&
$64$ & $6\cdot 192$\\
$(1,0,0,0)$ & $(\frac{m-n}{2},\frac{m+n}{2},\frac{m-n}{2},\frac{m+n}{2})$ &
$(\gvct{u}_1,\gvct{u}_0),(\gvct{u}_2,\gvct{u}_0),(\gvct{u}_0,\gvct{u}_1)$&
$32$ & $12\cdot 192$\\
\hline
$(1,0,0,0)$ & $(m,n,n,0)$ & $(\gvct{u}_1,\gvct{u}_0),(\gvct{u}_2,\gvct{u}_0),
(\frac{1}{\sqrt{2}}(0,1,1,0),\frac{1}{\sqrt{2}}(0,1,1,0))$ &
$32$ & $12 \cdot 192$\\
$(1,0,0,0)$ & $(\frac{m}{2}-n,\frac{m}{2}+n,\frac{m}{2},\frac{m}{2})$ & $(\gvct{u}_1,\gvct{u}_0),(\gvct{u}_2,\gvct{u}_0)$ &
$16$ & $24 \cdot 192$\\
\hline
$(1,0,0,0)$ & general & $(\gvct{u}_1,\gvct{u}_0),(\gvct{u}_2,\gvct{u}_0)$ &
$16$ & $24\cdot 192$\\
\hline
$(0,1,1,1)$ & $(0,1,1,1)$ & $
(\frac{1}{2}(1,1,1,1),\frac{1}{2}(1,1,1,1)),
(\frac{1}{\sqrt{2}}(0,1,-1,0),\frac{1}{\sqrt{2}}(0,1,-1,0))$
& $24$ & $16\cdot 192$\\
$(0,1,1,1)$ & $(3,1,1,1)$ & $
(\frac{1}{2}(1,1,1,1),\frac{1}{2}(1,1,1,1))$
& $12$ & $32\cdot 192$\\
\hline
$(0,1,1,1)$ & $(m,n,n,n)$ & $(\frac{1}{2}(1,1,1,1),\frac{1}{2}(1,1,1,1))$ & $12$ & $32\cdot 192$\\
\hline
$(0,1,1,1)$ & $(m,n,0,0)$ & $
(\frac{1}{\sqrt{2}}(0,0,1,-1),\frac{1}{\sqrt{2}}(1,1,0,0))$
& $16$ & $24\cdot 192$\\
$(0,1,1,1)$ & $(\frac{m-n}{2},\frac{m+n}{2},\frac{m-n}{2},\frac{m+n}{2})$ &
$(\gvct{u}_0,\gvct{u}_1)$ & $8$ & $48\cdot 192$\\
\hline
$(0,1,1,1)$ & $(m,n,n,0)$ & $
(\frac{1}{\sqrt{2}}(0,1,-1,0),\frac{1}{\sqrt{2}}(0,1,1,0))$
& $8$ & $48\cdot 192$\\
$(0,1,1,1)$ & $(\frac{m}{2}-n,\frac{m}{2}+n,\frac{m}{2},\frac{m}{2})$ &
--- & $4$ & $96\cdot 192$\\
\hline
$(0,1,1,1)$ & general & --- & $4$ & $96\cdot 192$\\
\hline
$(m,n,n,n)$ & $(m',n',n',n')$ & $(\frac{1}{2}(1,1,1,1),\frac{1}{2}(1,1,1,1))$ &
$12$ & $32\cdot 192$\\
\hline
$(m,n,n,n)$ & $(m',n',0,0)$ & $(\gvct{u}_0,\gvct{u}_1)$ & $8$ & $48\cdot 192$\\
\hline
$(m,n,n,n)$ & $(m',n',n',0)$ & --- & $4$ & $96\cdot 192$\\
\hline
$(m,n,n,n)$ & general & --- & $4$ & $96\cdot 192$\\
\hline
%
$(m,n,0,0)$ & $(m',n',0,0)$ &
$(\frac{1}{\sqrt{2}}(1,1,0,0),\frac{1}{\sqrt{2}}(1,1,0,0)),
(\frac{1}{\sqrt{2}}(1,1,0,0),\frac{1}{\sqrt{2}}(1,-1,0,0))$
& $32$ & $12\cdot 192$\\
$(m,n,0,0)$ & $(m',0,n',0)$ & $(\gvct{u}_1,\gvct{u}_0),(\gvct{u}_0,\gvct{u}_2)$
& $16$ & $24\cdot 192$\\
$(m,n,0,0)$ & $(\frac{m'-n'}{2},\frac{m'+n'}{2},\frac{m'-n'}{2},\frac{m'+n'}{2})$
& $(\gvct{u}_1,\gvct{u}_0),(\gvct{u}_0,\gvct{u}_1)$ & $16$ & $24\cdot 192$\\
$(m,n,0,0)$ & $(\frac{m'-n'}{2},\frac{m'+n'}{2},\frac{m'+n'}{2},\frac{m'-n'}{2})$
& $(\gvct{u}_1,\gvct{u}_0),(\gvct{u}_0,\gvct{u}_2)$ & $16$ & $24\cdot 192$\\
$(m,n,0,0)$ & $(\frac{m'-n'}{2},\frac{m'-n'}{2},\frac{m'+n'}{2},\frac{m'+n'}{2})$
& $(\gvct{u}_1,\gvct{u}_0),(\gvct{u}_0,\gvct{u}_3)$ & $16$ & $24\cdot 192$\\
\hline
$(m,n,0,0)$ & $(m',n',n',0)$ & $(\gvct{u}_1,\gvct{u}_0)$
& $8$ & $48\cdot 192$\\
$(m,n,0,0)$ & $(m',0,n',n')$ &
$(\frac{1}{\sqrt{2}}(1,1,0,0),\frac{1}{\sqrt{2}}(0,0,1,1))$
& $16$ & $24\cdot 192$\\
$(m,n,0,0)$ & $(\frac{m'}{2}-n',\frac{m'}{2}+n',\frac{m'}{2},\frac{m'}{2})$ &
$(\gvct{u}_1,\gvct{u}_0)$
& $8$ & $48\cdot 192$\\
$(m,n,0,0)$ & $(-\frac{m'}{2}-n',\frac{m'}{2},\frac{m'}{2}-n',\frac{m'}{2})$ & $(\gvct{u}_1,\gvct{u}_0)$
& $8$ & $48\cdot 192$\\
$(m,n,0,0)$ & $(\frac{m'}{2}-n',\frac{m'}{2},\frac{m'}{2}+n',\frac{m'}{2})$ & $(\gvct{u}_1,\gvct{u}_0)$
& $8$ & $48\cdot 192$\\
\hline
$(m,n,0,0)$ & general & $(\gvct{u}_1,\gvct{u}_0)$ & $8$ & $48\cdot 192$\\
\hline
$(m,n,n,0)$ & $(m',n',n',0)$ &
$(\frac{1}{\sqrt{2}}(0,1,1,0),\frac{1}{\sqrt{2}}(0,1,1,0))$
& $8$ & $48\cdot 192$\\
$(m,n,n,0)$ & $(m',0,n',n')$ & --- & $4$ & $96\cdot 192$\\
$(m,n,n,0)$ & $(\frac{m'}{2}-n',\frac{m'}{2}+n',\frac{m'}{2},\frac{m'}{2})$
& --- & $4$ & $96\cdot 192$\\
$(m,n,n,0)$ & $(-\frac{m'}{2}-n',\frac{m'}{2},\frac{m'}{2}-n',\frac{m'}{2})$
& --- & $4$ & $96\cdot 192$\\
$(m,n,n,0)$ & $(\frac{m'}{2}-n',\frac{m'}{2},\frac{m'}{2}+n',\frac{m'}{2})$
& --- & $4$ & $96\cdot 192$\\
\hline
$(m,n,n,0)$ & general & --- & $4$ & $96\cdot 192$\\
\hline
general & general & --- & $4$ & $96\cdot 192$\\
\hline
\end{tabular}
\end{sideways}
\caption{\label{tabebpartone} {P--E}quivalence classes of admissible pairs.
The groups $\CH_P(\gvct{p},\gvct{q})$ are given in terms of their
  generators listed in the third column, where we always have to add the
  trivial generators $((-1,0,0,0),(1,0,0,0))$ and $((1,0,0,0),(-1,0,0,0))$.
The last column gives the number of {P--e}quivalent coincidence
rotations $R(\gvct{p},\gvct{q})$. By dividing these numbers by 192 we obtain
the number of {P--e}quivalent CSLs.
In order to save space we have omitted some pairs. These
can be easily obtained by interchanging $\gvct{p}$ and
$\gvct{q}$ and adapting the generators correspondingly.}
\end{center}
\end{table}

\begin{table}[htb]
\begin{center}
\begin{tabular}{|c|c|c|}
\hline
$(\gvct{p},\gvct{q})$ &
\begin{minipage}{3cm}
$\#$ inequiv. pairs $\Sigma_P=\Sigma_F$
\end{minipage} &
\begin{minipage}{3cm}
$\#$ inequiv. pairs $\Sigma_P=2\Sigma_F$
\end{minipage}\\
\hline
$((1,0,0,0),(1,0,0,0))$& 1& ---\\
$(1,0,0,0),(1,1,1,1)$ & --- & 1\\
\hline
$((1,0,0,0),(m,n,n,n))$ & $n_{F02}$ & $2n_{F02}$\\
\hline
$((1,0,0,0),(m,n,0,0))$ & $n_{F03}$ & ---\\
$((1,0,0,0),(m-n,m+n,m-n,m+n))$ & --- & $n_{F03}$ \\
\hline
$((1,0,0,0),(m,n,n,0))$ & $n_{F04}$ & ---\\
$((1,0,0,0),(m-2n,m+2n,m,m))$ & --- & $n_{F04}$\\
\hline
$((1,0,0,0),(m,n,p,q))$ &$n_{F05}$ & $2n_{F05}$\\
\hline
$((0,1,1,1),(0,1,1,1))$ & 1 & --- \\
$((0,1,1,1),(3,1,1,1))$ & --- & 1 \\
\hline
$((0,1,1,1),(m,n,n,n))$ & $n_{F12}$ & $2n_{F12}$\\
\hline
$((0,1,1,1),(m,n,n,0))$ & $n_{F14}$ & --- \\
$((0,1,1,1),(m-2n,m+2n,m,m))$ & --- & $n_{F14}$\\
\hline
$((0,1,1,1) ,(m,n,p,q))$ & $n_{F15}$ & $2n_{F15}$\\
\hline
$((m,n,n,n),(m',n',n',n'))$ & $n_{F22}$ & $2n_{F22}$\\
\hline
$((m,n,n,n),(m',n',0,0))$ & $n_{F23}$ & $2n_{F23}$\\
\hline
$((m,n,n,n),(m',n',n',0))$ & $n_{F24}$ & $2n_{F24}$ \\
\hline
$((m,n,n,n),(m',n',p',q'))$ & $n_{F25}$ & $2n_{F25}$\\
\hline
$((m,n,0,0),(m',n',0,0))$ & $n_{F33}$ & ---\\
$((m,n,0,0),(m',0,n',0))$ & $n_{F33}$ & ---\\
$((m,n,0,0),(m'-n',m'+n',m'-n',m'+n'))$ & --- & $n_{F33}$ \\
$((m,n,0,0),(m'-n',m'+n',m'+n',m'-n'))$ & --- & $n_{F33}$ \\
$((m,n,0,0),(m'-n',m'-n',m'+n',m'+n'))$ & --- & $n_{F33}$ \\
\hline
$((m,n,0,0),(m',n',n',0))$ & $n_{F34}$ & ---\\
$((m,n,0,0),(m',0,n',n'))$ & $n_{F34}$ & ---\\
$((m,n,0,0),(m'-2n',m'+2n',m',m'))$ & --- & $n_{F34}$ \\
$((m,n,0,0),(-m'-2n',m',m'-2n',m'))$ & --- & $n_{F34}$ \\
$((m,n,0,0),(m'-2n',m',m'+2n',m'))$ & --- & $n_{F34}$ \\
\hline
$((m,n,0,0),(m',n',p',q'))$ & $3n_{F35}$ & $6n_{F35}$ \\
\hline
$((m,n,n,0),(m',n',n',0))$ & $n_{F44}$ &  ---\\
$((m,n,n,0),(m',0,n',n'))$ & $n_{F44}$ &  ---\\
$((m,n,n,0),(m'-2n',m'+2n',2,m'))$ & --- & $n_{F44}$ \\
$((m,n,n,0),(-m'-2n',m',m'-2n',m'))$ & --- &  $n_{F44}$ \\
$((m,n,n,0),(m'-2n',m',m'+2n',m'))$ & --- & $n_{F44}$ \\
\hline
$((m,n,n,0),(m',n',p',q'))$  & $3n_{F45}$ & $6n_{F45}$ \\
\hline
$((m,n,p,q) ,(m',n',p',q'))$ & $3n_{F55}$ & $6n_{F55}$ \\
\hline
\end{tabular}
\caption{\label{tabip} Number of \mbox{P--i}nequivalent admissible pairs.
The second column gives the number of inequivalent pairs for odd values of
$\Sigma$ whereas the third column gives the same information for even values
of $\Sigma$. In order to save space we have omitted some pairs. These
can be easily obtained by interchanging $\gvct{p}$ and $\gvct{q}$.}
\end{center}
\end{table}


\begin{thebibliography}{1}

\bibitem{boll70}
Bollmann, W.:
Crystal Defects and Crystalline Interfaces.
Springer, Berlin, 1970.

\bibitem{boll82}
Bollmann, W.:
Crystal lattices, interfaces, matrices.
published by the author, Geneva, 1982.

\bibitem{dun1}
Duneau, M.; Katz, A.:
Quasiperiodic patterns.
Phys. Rev. Lett. \textbf{54} (1985) 2688--2691.

\bibitem{mbaake98}
Baake, M.:
A Guide to Mathematical Quasicrystals
In: {\em Quasicrystals}, (Eds. J.-B.~Suck, M.~Schreiber, P.~H\"au{\ss}ler),
p. 17--48, Springer, Berlin, 2002.

\bibitem{baa97}
Baake, M.:
Solution of the coincidence problem in dimensions $d\leq 4$.
In: {\em {T}he {M}athematics of {L}ong-{R}ange
  Aperiodic Order} (Ed. R.~V. Moody), p. 9--44, Kluwer, Dordrecht, 1997.

\bibitem{pzcsl1}
Zeiner, P.:
Symmetries of coincidence site lattices of cubic lattices.
Z. Kristallogr. \textbf{220} (2005) 915--925.

\bibitem{koecheng}
Koecher, M.; Remmert, R.:
Hamilton's Quaternions.
In: {\em Numbers} (Eds. H.-D.~Ebbinghaus et. al.), p.~189--220.
  Springer, 1991.

\bibitem{hurw}
Hurwitz, A.:
Vor\-le\-sun\-gen \"uber die Zah\-len\-theo\-rie der Qua\-ter\-nio\-nen.
Springer, Berlin 1919.

\bibitem{val}
du~Val, P
Homographies, Quaternions and rotations.
Clarendon Press, Oxford, 1964.

\bibitem{grim74}
Grimmer, H.:
Disorientations and coincidence rotations for cubic lattices.
Acta Cryst. A {\bf 30} (1974) 685--688.

\bibitem{grim76}
Grimmer, H:
Coincidence site lattices: New results and comments on papers by
  Fortnow and Woirgard-de Fouquet.
Scripta Met. {\bf 10} (1976) 387--391.

\end{thebibliography}
\end{document}